\documentclass[12pt]{svmult}
\usepackage{amscd,amsfonts,amssymb,amsmath,latexsym,array,hhline}
\usepackage[dvips]{graphics}
\makeatletter 
\renewcommand{\@begintheorem}[2]{\begin{trivlist}\it
\item[\hspace{\labelsep}{\bf #1\ #2.}]}
\renewcommand{\@opargbegintheorem}[3]{\begin{trivlist}\it
\item[\hspace{\labelsep}{\bf #1\ #2\ (#3).}]}
\renewcommand{\@endtheorem}{\end{trivlist}}
\renewcommand{\@cite}[2]{[{#1\if@tempswa ; #2\fi}]}
\makeatother

\newtheorem{theo}{Theorem}
\newtheorem{prop}[theo]{Proposition}
\newtheorem{lemm}[theo]{Lemma}
\newtheorem{coro}[theo]{Corollary}

\newcommand{\paragr}{}

\renewcommand{\abstractname}{}

\newcommand\fdem{$\Box$}

\DeclareMathOperator{\nr}{\mathbb{N}}

\DeclareMathOperator{\re}{\mathbb{R}}
\DeclareMathOperator{\p}{\mathbb{P}}
\DeclareMathOperator{\e}{\mathbb{E}}

\DeclareMathOperator{\desp}{(\Omega, \mathcal{F}, (\mathcal F_n),\mathnormal{P})}

\author{{Jan Ob\l \'oj}\inst{1,2}\and{Marc Yor}\inst{1}\\}
\institute{
Laboratoire de Probabilit\'es et Mod\`eles Al\'eatoires, Universit\'e Paris 6,
\newline 4 pl. Jussieu, Bo\^{i}te 188, 75252 Paris Cedex 05, France
\and Wydzia\l\ Matematyki, Uniwersytet Warszawski
\newline Banacha 2, 02-097 Warszawa, Poland
\newline email: \texttt{obloj@mimuw.edu.pl}
}

\begin{document}
\title*{On local martingale and its supremum: harmonic functions and beyond.\\ {\small Dedicated to Professor A.N. Shiryaev for his 70$^{th}$ birthday}}
\titlerunning{On Brownian motion and its supremum}

\maketitle

\begin{abstract}
\textit{Version du 8 D\'ecembre 2004}\medskip\\
We discuss certain facts involving a continuous local martingale $N$ and its supremum $\overline{N}$. A complete characterization of $(N,\overline{N})$-harmonic functions is proposed. This yields an important family of martingales, the usefulness of which is demonstrated, by means of examples involving the Skorokhod embedding problem, bounds on the law of the supremum, or the local time at $0$, of a martingale with a fixed terminal distribution, or yet in some Brownian penalization problems. In particular we obtain new bounds on the law of the local time at $0$, which involve the excess wealth order.

.
\end{abstract}

\vfill
\newpage
\noindent\textbf{Dedication.} 
The first time I met Prof. A. Shiryaev was in January 1977, during a meeting dedicated to Control and Filtering theories, in Bonn. This was a time when meeting a Soviet mathematician was some event! Among the participants to that meeting, were, apart from A. Shiryaev, Prof. B. Grigelionis, and M. Yershov, who was by then just leaving Soviet Union in hard circumstances. To this day, I vividly remember that A.S, M.Y. and myself spent a full Sunday together, trying to solve a succession of problems raised by A.S., who among other things, explained at length about Tsirel'son's example of a one-dimensional SDE, with path dependent drift, and no strong solution (\cite{MR51:11654}; this motivated me to give - in \cite{MR91e:60175} - a more direct proof than the original one by Tsirel'son, see also \cite{MR93d:60104}, and Revuz and Yor \cite{MR2000h:60050} p. 392).

Each of my encounters with A.S. has had, roughly, the same flavor: A.S. would present, with great enthusiasm, some recent or not so recent result, and ask me for some simple proof, extension, etc... I have often been hooked into that game, which kept reminding me of one of my favorite pedagogical sentences by J. Dixmier: \textit{When looking for the 50$^{th}$ time at a well-known proof of some theorem, I would discover a new twist I had never thought of, which would cast a new light on the matter}.

I hope that the following notes, which discuss some facts about local martingales and their supremum processes, and are closely related to the thesis subject of the first author, may also have some this ``new twist'' character for some readers, and be enjoyed by Albert Shiryaev, on the occasion of his 70$^{th}$ birthday.
\begin{flushright}
  \textit{Marc Yor}
\end{flushright}
\section{Introduction}
\label{sec:intro}

In this article we focus on local martingales, functions of two-dimensional processes, whose components are a continuous local martingale $(N_t:t\geq 0)$ and its supremum $\overline{N}_t=\sup_{u\le t}N_u$, i.e. on local martingales of the form $(H(N_t,\overline{N}_t):t\ge 0)$, where  $H:\re\times\re_+\to\re$. We call functions $H:\re\times\re_+\to\re$ such that $(H(N_t,\overline{N}_t):t\ge 0)$ is a local martingale, $(N,\overline{N})$-harmonic functions.

Some examples of such local martingales are 
\begin{equation}
  \label{eq:Ffmart}
F(\overline{N}_t)-f(\overline{N}_t)(\overline{N}_t-N_t),\quad t\ge 0,  
\end{equation}
where $F\in C^1$ and $F'=f$, introduced by Az\'ema and Yor \cite{MR82c:60073a}. 
We show that (\ref{eq:Ffmart}) defines a local martingale for any Borel, locally integrable function $f$.
We conjecture that these are the only local martingales, that is that the only $(N,\overline{N})$-harmonic functions are of the form $H(x,y)=F(y)-f(y)(y-x)+C$, with $f$ a locally integrable function, $F(y)=\int_0^y f(u)du$, and $C$ a constant.

We explain, in an intuitive manner, how these local martingales, which we call \emph{max-martingales}, may be used to find the Az\'ema-Yor solution to the Skorokhod embedding problem. We then go on and develop, with the help of  these martingales, the well-known bounds on the law of the supremum of a uniformly integrable martingale with a fixed terminal distribution. Using the L\'evy and Dambis-Dubins-Schwarz theorems, we reformulate the results in terms of the absolute value $|N|$ and the local time $L^N$ at $0$, of the local martingale $N$. This leads to some new bounds on the law of the local time of a uniformly integrable martingale with fixed terminal distribution. A recently introduced and studied stochastic order, called the excess wealth order (see Shaked and Shanthikumar \cite{MR98j:60028}), plays a crucial role. We also point out that the max-martingales appear naturally in some Brownian penalization problems. Finally, we try to sketch a somewhat more general viewpoint linked with the balayage formula.

The organization of this paper is as follows.
We start in Section \ref{sec:discrete} with a discrete version of the balayage formula and show how to deduce from it Doob's maximal and $L^p$ inequalities. In the subsequent Section \ref{sec:harm}, in Theorem \ref{thm:Bharm}, we formulate the result about the harmonic functions of $(N,\overline{N})$ and prove it in a regular case. Section \ref{sec:appl} is devoted to some applications: it contains three subsections concentrating respectively on the Skorokhod embedding problem, bounds on the laws of $\overline{N}$ and $L^N$, and Brownian penalizations.
The last section contains a discussion of the balayage formula.

\section{Notation}
\label{sec:notation}

Throughout this paper $(N_t:t\geq 0)$ will denote a continuous local martingale, and $\overline{N}_t=\sup_{s\le t}N_s$. We have adopted this notation so that there is no confusion with stock-price processes, which are often denoted $S_t$. The local time at $0$ of $N$ is denoted $(L^N_t:t\ge 0)$. For processes either in discrete or in continuous time, when we say that a process is a (sub/super) martingale without specifying the filtration, we mean the natural filtration of the process.

$B=(B_t:t\ge 0)$ shall denote a one-dimensional Brownian motion, starting from $0$, and $\overline{B}_t=\sup_{s\leq t}B_s$. The natural filtration of $B$ is denoted $(\mathcal{F}_t)$ and is taken completed.

The indicator function is denoted $\mathbf{1}$. We define $\lor$ and $\land$ through $a\lor b=\max\{a,b\}$ and $a\land b=\min\{a,b\}$. The positive part is given by $x^+=x\lor 0$.  
 For $\mu$ a probability measure on $\re$, $\overline{\mu}(x):=\mu([x,\infty))$ is its tail distribution function. $X\sim \mu$ means \textit{X has distribution $\mu$.}

\section{Balayage in discrete time and some applications}
\label{sec:discrete}

We start with the discrete time setting, and present a simple idea, which corresponds to balayage in continuous time, and which proves an efficient tool, as it allows, for example, to obtain easily Doob's maximal and $L^p$ inequalities. Let $\desp$ be a filtered probability space and $(Y_n:n\geq 0)$ be some real-valued adapted discrete stochastic process. Let $(\phi_n: n\geq 0)$ be also an adapted process, which further satisfies $\phi_n\mathbf{1}_{Y_n\neq 0}=\phi_{n-1}\mathbf{1}_{Y_n\neq 0}$, for all $n\in \nr$. The last condition can be also formulated as ``the process $(\phi_n)$ is constant on excursions of $(Y_n)$ away from $0$''.
\begin{lemm}
\label{lem:dbal}
  Let $(Y_n,\phi_n)$ be as above, $Y_0=0$. The following identities hold:
  \begin{equation}
    \label{eq:d_bal}
    \phi_n Y_n= \phi_{n-1} Y_n=\sum_{k=1}^n \phi_{k-1}(Y_k-Y_{k-1}),\quad n\geq 1.
  \end{equation}
\end{lemm}
{\sl Proof.}
  The first equality is obvious as $\phi_n Y_n=\phi_n Y_n\mathbf{1}_{Y_n \neq 0}=\phi_{n-1} Y_n$, and the second one is obtained by telescoping.
\fdem

To see how the above can be used, let us give some examples of pairs $(Y_n,\phi_n)$ involving in particular an adapted process $X_n$ and its maximum $\overline{X}_n$:
\begin{itemize}
\item $Y_n=\overline{X}_n-X_n$ and $\phi_n=f(\overline{X}_n)$, for some Borel function $f$;
\item $Y_n=X_n$, $\phi_n=\sum_{k=0}^n \mathbf{1}_{X_k=0}$ (note that $Y_n=|X_n|$ works as well);
\item $Y_n=X^*_n-|X_n|$, $\phi_n=f(X^*_n)$, for some Borel function $f$, where $X^*_n=\max_{k\le n}|X_k|$;
\item $Y_n=\overline{X}_n-\underline{X}_n$, $\phi_n=f\Big(\sum_{k=1}^n \mathbf{1}_{(\overline{X}_k=\underline{X}_k)}\Big)$, for some Borel function $f$, where $\underline{X}_n=|\min_{k\le n} X_k|$.
\end{itemize}
We now use the discrete balayage formula with the first of the above examples to establish a useful supermartingale property.
\begin{prop}
\label{prop:discmart}
  Let $(X_n:n\in\nr)$ be a submartingale in its natural filtration $(\mathcal{F}_n)$, $X_0=0$, and $f$ some increasing, locally integrable, positive function. 
Assume that $\e f(\overline{X_n})<\infty$ and $\e F(\overline{X}_n)<\infty$ for all $n\in\nr$,  where $F(x)=\int_0^x f(s)ds$.
Then the process $S^f_n=f(\overline{X}_n)(\overline{X}_n-X_n)-F(\overline{X}_n)$ is a $(\mathcal{F}_n)$-supermartingale.
\end{prop}
{\sl Proof.}
  Since the pair $(\overline{X}_n-X_n, f(\overline{X}_n))$ satisfies the assumptions of Lemma \ref{lem:dbal}, we have:
  \begin{eqnarray}
\label{eq:disc_pom}    
S^f_n&=&\sum_{k=1}^nf(\overline{X}_{k-1})(\overline{X}_k-X_{k}-\overline{X}_{k-1}+X_{k-1})-F(\overline{X}_n)\nonumber\\
    &=&\sum_{k=1}f(\overline{X}_{k-1})(\overline{X}_{k}-\overline{X}_{k-1})-\sum_{k=1}^n f(\overline{X}_{k-1})(X_k-X_{k-1})-\int_0^{\overline{X}_n}f(x)dx\nonumber\\
    &=&\sum_{k=1}^n\int_{\overline{X}_{k-1}}^{\overline{X}_k}\Big(f(\overline{X}_{k-1})-f(x)\Big)dx - \sum_{k=1}^n f(\overline{X}_{k-1})(X_k-X_{k-1}).
  \end{eqnarray}
Using (\ref{eq:disc_pom}), the fact that $f$ is increasing, and $(X_n)$ is a submartingale, we obtain the supermartingale property for $S^f_n$.
\fdem

The above Proposition allows to recover Doob's maximal and $L^p$ inequalities in a very easy way. Indeed, consider the function $f(x)=\mathbf{1}_{x\ge \lambda}$ for some $\lambda>0$. Then the process $S^f_n=S^{(\lambda)}_n=\mathbf{1}_{\overline{X}_n\ge \lambda}(\lambda-X_n)$ is a supermartingale, which yields \textbf{Doob's maximal inequality}
\begin{equation}
  \label{eq:doobmax}
  \lambda \p\Big(\overline{X}_n\ge \lambda\Big)\le \e \Big[\mathbf{1}_{(\overline{X}_n\ge \lambda)}X_n\Big].
\end{equation}

To obtain the $L^p$ inequalities  we consider the function $f(x)=px^{p-1}$ for some $p>1$, and we suppose that $(X_n:n\geq 0)$ is a positive submartingale with $\e X_n^p<\infty$. This implies, as $\overline{X}_n^p\le \sum_{k=1}^nX_k^p$, that $\e \overline{X}_n^p<\infty$.
The process $S^f_n=S^{(p)}_n=(p-1)(\overline{X}_n)^p-p(\overline{X}_n)^{p-1}X_n$ is then a supermartingale, which yields
\begin{eqnarray}
  \label{eq:doobslp}
  (p-1)\e\Big[(\overline{X}_n)^p\Big]&\le& p \e\Big[(\overline{X}_{n})^{p-1}X_n\Big]\quad\textrm{and hence, applying H\"older's inequality,}\nonumber\\
\e\Big[(\overline{X}_n)^p\Big]&\leq& \Big(\frac{p}{p-1}\Big)^p\e\Big[X_n^p\Big],\quad \textrm{ which is }\textbf{Doob's }\mathbf{L^p}\textbf{ inequality}.
\end{eqnarray}

To our best knowledge, this small wrinkle about Doob's inequalities for positive submartingales involving supermartingales does not appear in any of the books on discrete martingales, such as Neveu \cite{MR53:6729}, Garsia \cite{MR56:6844} or Williams \cite{MR93d:60002}. We point out also, that our method allows to obtain other variants of Doob's inequalities, such as $L\log L$ inequalities, etc.

\section{The Markov process $((B_t,\overline{B}_t):t\ge 0)$ and its harmonic functions}
\label{sec:harm}

In the rest of the paper we will focus on the continuous-time setup. It follows immediately from the strong Markov property of $B$, or rather the independence of its increments, that for $s<t$, and $f:\re\times\re_+\to\re_+$ a Borel function, one has:
\begin{equation}
  \label{eq:markovbb}
  \e\Big[f(B_t,\overline{B}_t)\Big|\mathcal{F}_s\Big]=\tilde{\e}\Big[f\Big(B_s+\tilde{B}_{t-s},\overline{B}_s\lor \sup_{u\le t-s}(B_s+\tilde{B}_u)\Big)\Big],
\end{equation}
where on the RHS, the notation $\tilde{\e}$ indicates integration with respect to functionals of the Brownian motion $(\tilde{B}_u:u\ge 0)$, which is assumed to be independent of $(B_t:t\ge 0)$.

In particular, the two-dimensional process $((B_t,\overline{B}_t): t\ge 0)$ is a nice Markov process, hence a strong Markov process, and its semigroup can be computed explicitly thanks to the well-known, and classical formula:
\begin{equation*}
  \label{eq:lawbb}
  \p\Big(B_t\in dx,\overline{B}_t\in dy\Big)=\Big(\frac{2}{\pi t^3}\Big)^{1/2}(2y-x)\exp\Big(-\frac{(2y-x)^2}{2t}\Big)\mathbf{1}_{(y\ge x^+)}dxdy.
\end{equation*}

We are now interested in a description of the harmonic functions $H: \nolinebreak[4]\re\times \re_+\to \re$ of $(B,\overline{B})$ that is of Borel functions such that $(H(B_t,\overline{B}_t):t\ge 0)$ is a local martingale. Note that this question is rather natural and interesting since $H$ is $(B,\overline{B})$-harmonic if and only if, thanks to the Dambis-Dubins-Schwarz theorem, for any continuous local martingale $(N_t:t\ge 0)$, $H$ is also $(N,\overline{N})$-harmonic. The following Proposition is an extension of Proposition 4.7 in Revuz and Yor \cite{MR2000h:60050}.

\begin{theo}
\label{thm:Bharm}
  Let $N=(N_t:t\ge 0)$ be a continuous local martingale with $\langle N\rangle_\infty=\infty$ a.s., $f$ a Borel, locally integrable function, and $H$ defined through
  \begin{equation}
    \label{eq:Hgen}
    H(x,y)=F(y)-f(y)(y-x)+C,
  \end{equation}
where $C$ is a constant and $F(y)=\int_0^y f(s)ds$. Then, the following holds:
\begin{equation}
  \label{eq:Hgenmart}
H(N_t,\overline{N}_t)=F(\overline{N}_t)-f(\overline{N}_t)(\overline{N}_t-N_t)+C=\int_0^t f(\overline{N}_s)dN_s+C, \quad t\geq 0,  
\end{equation}
and $(H(N_t,\overline{N}_t):t\ge 0)$ is a local martingale.
\end{theo}
\textbf{Remarks.} Local martingales of the form (\ref{eq:Hgenmart}) were first introduced by Az\'ema and Yor \cite{MR82c:60073a} and used to solve the Skorokhod embedding problem (cf. Section \ref{sec:skoro} below). In the light of the above theorem, we will call them \textbf{max-martingales} and the functions given in (\ref{eq:Hgen}) will be called \textbf{MM-harmonic functions} (max-martingale harmonic) or $(N,\overline{N})$-harmonic. Note the resemblance of (\ref{eq:Hgenmart}) with the discrete time process $S^f_n$ given in Proposition \ref{prop:discmart}.\\
It is known (see Proposition 4.7 in Revuz and Yor \cite{MR2000h:60050}) that if $H\in C^{2,1}$ then the reverse statement holds. That is, if $H$ is $(N,\overline{N})$-harmonic then there exists a continuous function $f$ such that (\ref{eq:Hgen}) holds. We present below a proof of this fact. We conjecture that it suffices to suppose only that $H$ is right-continuous in the first coordinate. We plan to develop this in a separate work.
\smallskip\\
{\sl Proof.}
As mentioned above, thanks to the Dambis-Dubins-Schwarz theorem, it suffices to prove the theorem for $N=B$.  
We first recall how to prove the converse of the theorem for the regular case. We assume that $H\in C^{2,1}$, with obvious notation, and that $H$ is $(B,\overline{B})$-harmonic.
We denote by $H_x'$ and $H_y'$ the partial derivatives of $H$ in the first and the second argument respectively, and $H_{x^2}''$ the second derivative of $H$ in the first argument.
Without loss of generality, we assume that $H(0,0)=0$.
Under the present assumptions we can apply It\^o's formula to obtain:
\begin{equation*}
  \label{eq:itoH}
  H(B_t,\overline{B}_t)=\int_0^t H_x'(B_s,\overline{B}_s)dB_s+\int_0^tH_y'(\overline{B}_s,\overline{B}_s)d\overline{B}_s+\frac{1}{2}\int_0^t H_{x^2}''(B_s,\overline{B}_s)ds,
\end{equation*}
where we used the fact that $B_s=\overline{B}_s$, $d\overline{B}_s$-a.s. Now, since $H(B_t,\overline{B_t})$ is a local martingale, the above identity holds if and only if:
\begin{equation}
  \label{eq:condH}
  H_y'(\overline{B}_s,\overline{B}_s)d\overline{B}_s+\frac{1}{2}H_{x^2}''(B_s,\overline{B}_s)ds=0,\quad s\geq 0.
\end{equation}
The random measures $d\overline{B}_s$ and $ds$ are mutually singular since $d\overline{B}_s=\mathbf{1}_{(\overline{B}_s-B_s=0)}d\overline{B}_s$ and $ds=d\langle B\rangle_s=\mathbf{1}_{(\overline{B}_s-B_s\neq 0)}d\langle B\rangle_s$. Equation (\ref{eq:condH}) holds therefore if and only if
\begin{equation}
  \label{eq:condH2}
  H_y'(y,y)=0\quad\textrm{and}\quad H_{x^2}''(x,y)=0.
\end{equation}
The second condition implies that $H(x,y)=f(y)x+g(y)$ and the first one then gives $f'(y)y+g'(y)=0$. Thus, $g(y)=-\int_0^y uf'(u)du=\int_0^y f(u)du-f(y)y$. This yields formula (\ref{eq:Hgen}).

Furthermore, the above reasoning grants us that the formula (\ref{eq:Hgenmart}) holds for $f$ of class $C^1$. As $C^1$ is dense in the class of locally integrable functions (in an appropriate norm), if we can show that 
the quantities given in (\ref{eq:Hgenmart}) are well defined and finite for any locally integrable $f$ on $[0,\infty)$, then the formula (\ref{eq:Hgenmart}) extends to such functions through monotone class theorems. 
For $f$ a locally integrable function, $F(x)$ is well defined and finite, so all we need to show is that $\int_0^t f(\overline{B}_s)dB_s$ is well defined and finite a.s. This is equivalent to $\int_0^t \Big(f(\overline{B}_s) \Big)^2ds<\infty$ a.s., which we now show.

Write $T_x=\inf\{t\ge 0: B_t=x\}$ for the first hitting time of $x$, which is a well defined, a.s. finite, stopping time. Thus $\int_0^t \Big(f(\overline{B}_s) \Big)^2ds<\infty$ a.s., if and only if, for all $x>0$, $\int_0^{T_x} \Big(f(\overline{B}_s) \Big)^2ds<\infty$. However, the last integral can be rewritten as
\begin{eqnarray}
  \label{eq:localvsqlocal}
  \int_0^{T_x}ds\Big(f(\overline{B}_s)\Big)^2&=&\sum_{0\le u\le x}\int_{T_{u-}}^{T_u}ds \Big(f(\overline{B}_s)\Big)^2\nonumber\\
&=&\sum_{0\le u \le x}f^2(u)\Big(T_u-T_{u-}\Big)=\int_0^x f^2(u)d T_u.
\end{eqnarray}
Now it suffices to note that (see Ex. III.4.5 in Revuz and Yor \cite{MR2000h:60050})
\begin{equation}
  \e\Big[\exp\Big(-\frac{1}{2}\int_0^xf^2(u)dT_u\Big)\Big]=\exp\Big(-\int_0^x|f(u)|du\Big),
\end{equation}
to see that the last integral in (\ref{eq:localvsqlocal}) is finite if and only if $\int_0^x|f(u)|du<\infty$, which is precisely our hypothesis on $f$.

Note that the function $H$ given by (\ref{eq:Hgen}) is locally integrable as both $x\to f(x)$ and $x\to xf(x)$ are locally integrable. 
\fdem

L\'evy's theorem guarantees that the processes $((B_t,\overline{B}_t):t\geq 0)$ and $((L_t-|B_t|,L_t):t\geq 0)$ have the same distribution, where $L_t$ denotes local time at $0$ of $B$. Theorem \ref{thm:Bharm} yields therefore also a complete description of $(L,|B|)$-harmonic functions, which again through Dambis-Dubins-Schwarz theorem, extends to any local continuous martingale. We have the following
\begin{coro}
  \label{thm:Lharm}
Let $N=(N_t:t\ge 0)$ be a continuous local martingale with $\langle N\rangle_\infty=\infty$ a.s., and $L^N=(L_t^N:t\ge 0)$ its local time at $0$. Let  $g$ a Borel, locally integrable function, and $H$ be defined through
\begin{equation}
  \label{eq:HgenL}
  H(x,y)=G(y)-g(y)x+C,
\end{equation}
where $C$ is a constant and $G(y)=\int_0^y g(s)ds$. Then, the following holds:
\begin{equation}
  \label{eq:HgenmartL}
H(|N_t|,L^N_t)=G(L^N_t)-g(L^N_t)|N_t|+C=-\int_0^t g(L^N_s)sgn(N_s)dN_s+C, \quad t\geq 0,  
\end{equation}
and $(H(|N_t|,L^N_t):t\ge 0)$ is a local martingale.
\end{coro}

\section{Some appearances of the MM-harmonic functions}
\label{sec:appl}

We now present some easy applications of the martingales described in the previous section. We will show an intuitive way to obtain a solution to the Skorokhod embedding problem, as given by Az\'ema and Yor \cite{MR82c:60073a}. We will also discuss relations between the law of $\overline{B}_T$ and the conditional law of $B_T$ knowing $\overline{B}_T$, for some stopping time $T$. In the second subsection we will derive well-known bounds on the law of $\overline{B}_T$, when the law of $B_T$ is fixed. We will then continue in the same vein and describe the law of $L_T$, when the law of $B_T$ is fixed. We will end with a discussion of penalization of Brownian motion with a function of its supremum and some absolute continuity relations.

\subsection{On the Skorokhod embedding problem}
\label{sec:skoro}

The classical Skorokhod embedding problem can be formulated as follows: for a given centered probability measure $\mu$, find a stopping time $T$ such that $B_T\sim \mu$ and $(B_{t\land T}:t\geq 0)$ is a uniformly integrable martingale. Numerous solutions to this problem are known; for an extensive survey see Ob\l \'oj \cite{genealogia}. Here we make a remark about the solution given by Az\'ema and Yor in \cite{MR82c:60073a}. Namely we point out how one can arrive intuitively to this solution using the max-martingales (\ref{eq:Hgenmart}). Naturally, this might be extracted from the original paper, but it may not be so obvious to do so.

The basic observation is that the max-martingales allow to express the law of the terminal value of $\overline{B}$, that is $\overline{B}_T$, in terms of the conditional distribution of $B_T$ given $\overline{B}_T$. One then constructs a stopping time which actually binds both terminal values through a function and sees that the function can be obtained in terms of the target measure $\mu$.
\begin{prop}[Vallois \cite{MR95f:60051}]
\label{prop:lawofsup}
  Let $T$ be a $(\mathcal{F}_t)$-stopping time, such that $(B_{t\land T}:t\ge 0)$ is a uniformly integrable martingale. Write $\nu$ for the law of $\overline{B}_T$ and suppose that $\nu$ is equivalent to the Lebesgue measure on its interval support $[0,b]$, $b\leq \infty$.
Then the law of $\overline{B}_T$ is given by:
  \begin{equation}
    \label{eq:lawsupB}
    \p\Big(\overline{B}_T\geq y\Big)=\exp\Big(-\int_0^y\frac{ds}{s-\varphi(s)}\Big),\quad 0\leq y\leq b,
  \end{equation}
where $\varphi(x)=\e[B_T|\overline{B}_T=x]$, i.e. $\e[B_Th(\overline{B_T})]=\e[\varphi(\overline{B}_T)h(\overline{B}_T)]$, for any positive Borel function $h$.
\end{prop}
\textbf{Remark.} We note that the above formula in the special case when $B_T=\varphi(\overline{B}_T)$ a.s., and actually in the more general context of time-homogeneous diffusions, was obtained already by Lehoczky \cite{MR56:16770}. Vallois \cite{MR95f:60051} studied this issues in detail and has some more general formulae.
\\
{\sl Proof.}
  With the help of the max-martingales, we get, for any $f:\re_+\to\re$, bounded, with compact support:
  \begin{equation*}
    \e\Big[F(\overline{B}_T)-f(\overline{B}_T)(\overline{B}_T-B_T)\Big]=0.
  \end{equation*}
Upon conditioning with respect to $\overline{B}_T$ we obtain:
\begin{equation}
  \label{eq:pomoc1}
    \e\Big[F(\overline{B}_T)-f(\overline{B}_T)(\overline{B}_T-\varphi(\overline{B}_T))\Big]=0.
\end{equation}
We can rewrite the above as a differential equation involving $\nu\sim \overline{B}_T$, which yields (\ref{eq:lawsupB}).
\fdem

Define the Az\'ema-Yor stopping time, as suggested above, through $T_{\varphi}=\inf\{t\ge 0: B_t=\varphi(\overline{B}_t)\}$, for some strictly increasing, continuous function $\varphi:\re_+\to\re$. Obviously $B_{T_\varphi}=\varphi(\overline{B}_{T_\varphi})$. We look for a function $\varphi=\varphi_\mu$ such that $B_{T_{\varphi_\mu}}\sim \mu$. To this end, we take $x$ in the support of $\mu$ and write
\begin{equation*}
  \overline{\mu}(x)=\p(B_{T_{\varphi_\mu}}\ge x)=\overline{\nu}(\varphi_\mu^{-1}(x))=\exp\Big(-\int_0^{\varphi_\mu^{-1}(x)}\frac{ds}{s-\varphi_\mu(s)}\Big),
\end{equation*}
which may be considered as an equation on $\varphi_\mu$ in terms of $\mu$. Solving this equation, one obtains
\begin{equation}
  \label{eq:HL}
  \varphi^{-1}_\mu(x)=\Psi_\mu(x)=\frac{1}{\overline{\mu}(x)}\int_{[x,\infty)}sd\mu(s),
\end{equation}
the Hardy-Littlewood maximal function, or barycentre function, of $\mu$. 
\begin{prop}[Az\'ema-Yor \cite{MR82c:60073a}]
\label{prop:ay}
  Let $\mu$ be a centered probability measure. Define the function $\Psi_\mu$ through (\ref{eq:HL}) for $x$ such that $\overline{\mu}(x)\in (0,1)$ and put $\Psi_\mu(x)=0$ for $x$ such that $\overline{\mu}(x)=1$, $\Psi_\mu(x)=x$ for $x$ such that $\overline{\mu}(x)=0$. Then the stopping time $T_\mu:=\inf\{t\ge 0: \overline{B}_t=\Psi_\mu(B_t)\}$ satisfies $B_{T_\mu}\sim\mu$ and $(B_{t\land T_\mu}:t\ge 0)$ is a uniformly integrable martingale.
\end{prop}
The arguments presented above contain the principal ideas behind the Az\'ema-Yor solution to the Skorokhod embedding problem. Naturally, they work well for measures with positive density on $\re$. A complete proof of Proposition \ref{prop:ay} requires some rigorous arguments involving, for example, a limit procedure, but this can be done, as shown by Michel Pierre \cite{MR82f:60108}.

We now develop a link between formula (\ref{eq:lawsupB}) and work of Rogers \cite{MR94h:60067}. Let us carry out some formal computations. Write $\rho$ for the law of the couple $(\overline{B}_T, \overline{B}_T-B_T)\in \re_+\times\re_+$, and $\nu$ for its first marginal (as above). Differentiating (\ref{eq:lawsupB}) we  find
\begin{eqnarray}
  \label{eq:rogers}
 && d\overline{\nu}(y)\;\;=-\frac{\overline{\nu}(y)dy}{y-\varphi(y)},\quad\textrm{hence}\nonumber\\
&&\overline{\nu}(y)dy=(y-\varphi(y))d\nu(y),\quad\textrm{which we rewrite in terms of }\rho\nonumber\\
&&\Big(\iint_{(y,\infty)\times\re_+}\rho(ds,dx)\Big)dz=\int_{(0,\infty)}z\rho(ds,dz).
\end{eqnarray}
The last condition appears in Rogers \cite{MR94h:60067} and is shown to be equivalent to the existence of a continuous, uniformly integrable martingale $(B_{t\land T}:t\ge 0)$ such that $(\overline{B}_T, \overline{B}_T-B_T)\sim \rho$. Our formulation in (\ref{eq:lawsupB}) is less general, as it is not valid when the law of $\overline{B}_T$ has atoms. However, when it is valid, it provides an intuitive reading of (\ref{eq:rogers}).

To close this section, we point out that arguments similar to the ones presented above, can be developed to obtain a solution to the Skorokhod embedding problem for $|B|$ based on $L$: it suffices to use the martingales given by (\ref{eq:HgenmartL}) instead of those given by (\ref{eq:Hgenmart}). For a probability measure $m$ on $\re_+$, define the dual Hardy-Littlewood function (see Ob\l \'oj and Yor \cite{obloj_yor}) through
\begin{equation}
  \label{eq:dualHL}
  \psi_m(x)=\int_{[0,x]}\frac{y}{\overline{m}(y)}dm(y),\quad \textrm{for }x\textrm{ such that }\overline{m}(x)\in(0,1),
\end{equation}
and put $\psi_m(x)=0$ for $x$ such that $\overline{m}(x-)=1$ and $\psi_m(x)=\infty$ for $x$ such that $\overline{m}(x+)=0$.
\begin{prop}[Vallois \cite{MR86m:60200}, Ob\l \'oj and Yor \cite{obloj_yor}]
\label{prop:oy}
  Let $m$ be a probability measure on $\re_+$ with $m(\{0\})=0$ and define the function $\psi_m$ through (\ref{eq:dualHL}). Let $\varphi_m(y)=\inf\{x\geq 0: \psi_m>y\}$ be the right inverse of $\psi_m$.
Then the stopping time $T^m:=\inf\{t> 0: |B|_t=\varphi_m(L_t)\}$ satisfies $|B|_{T^m}\sim m$. Furthermore, $(B_{t\land T^m}:t\ge 0)$ is a uniformly integrable martingale if and only if $\int_0^\infty x dm(x)<\infty$.
\end{prop}
We note that the law of $L_{T^m}$ is given through $\p(L_{T^m}\geq x)=\exp\Big(-\int_0^x \frac{ds}{\varphi_m(s)}\Big)$ (cf. $(5.4)$ in \cite{obloj_yor}). An easy analogue of Proposition \ref{prop:lawofsup}, is that this formula is also true for general stopping time $T$, such that the law of $L_T$ has a density, with the function $\varphi_m$ replaced by $\varphi(x)=\e[|B_T||L_T=x]$.

\subsection{Bounds on the laws of $\overline{B}_T$ and $L_T$}
\label{sec:bounds}

We present a classical bound on the law of $\overline{B}_T$, which was first obtained by Blackwell and Dubins \cite{MR27:1557} and Dubins and Gilat \cite{MR58:13333} (see also Az\'ema and Yor \cite{MR82c:60073b}, Kertz and R\"osler \cite{MR91g:60056} and Hobson \cite{MR2000c:60056}).
\begin{prop}
Let $\mu$ be a centered probability measure and $T$ a stopping time, such that $B_T\sim \mu$ and $(B_{t\land T}:t\ge 0)$ is a uniformly integrable martingale. Then the following bound is true
\begin{equation}
  \label{eq:bound_sup}
  \p(\overline{B}_T\ge \lambda)\leq \p(\overline{B}_{T_\mu}\ge \lambda)=\overline{\mu}(\Psi_\mu^{-1}(\lambda)),\quad \lambda\geq 0
\end{equation}
where $T_\mu$ is given in Proposition \ref{prop:ay}, $\Psi_\mu$ is displayed in (\ref{eq:HL}) and its inverse is taken right-continuous.
\end{prop}
In other words, for the partial order given by tails domination, the law of $\overline{B}_T$ is bounded by the image of $\mu$ through the Hardy-Littlewood maximal function (\ref{eq:HL}).\\
{\sl Proof.}
 Suppose for simplicity that $\mu$ has a positive density, which is equivalent to $\Psi_\mu$ being continuous and strictly increasing. We consider the max-martingale (\ref{eq:Hgenmart}) for $f(x)=\mathbf{1}_{(x\ge \lambda)}$, for some fixed $\lambda>0$, and apply the optional stopping theorem. We obtain:
\begin{equation}
  \label{eq:doobmaxcont}
  \lambda P(\overline{B}_T\ge \lambda)= \e\Big[B_T\mathbf{1}_{(\overline{B}_T\ge \lambda)}\Big],
\end{equation}
that is Doob's maximal equality for continuous-time martingales. Let $p:=\p(\overline{B}_T\ge \lambda)$. As $B_T\sim \mu$, then the RHS is smaller than $\e[B_T \mathbf{1}_{(B_T\ge \overline{\mu}^{-1}(p))}]$ which, by definition in (\ref{eq:HL}), is equal to $p\Psi_\mu(\overline{\mu}^{-1}(p))$. We obtain therefore:
\begin{eqnarray}
  \lambda\p(\overline{B}_T\ge \lambda)=\lambda p&\leq& \e\Big[B_T\mathbf{1}_{(B_T\ge \overline{\mu}^{-1}(p))}\Big]=p\Psi_\mu\Big(\overline{\mu}^{-1}(p)\Big),\quad\textrm{hence}\nonumber\\
p&\leq& \overline{\mu}\Big(\Psi_\mu^{-1}(\lambda)\Big),\quad\textrm{since }\overline{\mu}\textrm{ is decreasing}.
\end{eqnarray}
To end the proof is suffice to note that $\p(\overline{B}_{T_\mu}\ge \lambda)=\overline{\mu}(\Psi_\mu^{-1}(\lambda))$, which is obvious from the definition of $T_\mu$. 
\fdem

Investigation of similar quantities with $\overline{B}_T$ replaced by $T$ is also possible. 
 Numerous authors studied the limit $\sqrt{\lambda}\p(T\ge \lambda)$.
It goes back to Az\'ema, Gundy and Yor \cite{MR580108} with more recent works by Elworthy, Li and Yor \cite{MR1478722} and Peskir and Shiryaev \cite{MR99j:60123}.

Integrating (\ref{eq:bound_sup}) one obtains bounds on the expectation of $\overline{B}_T$. Another bound on $\e \overline{B}_T$ can be obtained using the max-martingales. Take $f(x)=2x$, then by (\ref{eq:Hgenmart}) the process $B_t^2-2\overline{B}_tB_t=(\overline{B}_t-B_t)^2-B^2_t$ is a local martingale. For a bounded stopping time $T$, we have then $\e (\overline{B}_T-B_T)^2=\e B^2_T$, which yields:
\begin{equation}
\label{eq:max1}
  \e \overline{B}_T=\e(\overline{B}_T-B_T)\le \sqrt{\e\Big[(\overline{B}_T-B_T)^2\Big]}=\sqrt{\e\Big[B_T^2\Big]}=\sqrt{\e T}. 
\end{equation}
The inequality $\e \overline{B}_T\le \sqrt{\e T}$ extends to any stopping time, through the monotone convergence theorem. This inequality was generalized for Bessel processes by Dubins, Shepp and Shiryaev \cite{MR96j:60077} and for Brownian motion with drift by Peskir and Shiryaev \cite{MR2002i:60150}. These problems are also in close relation with the so-called Russian options developed mainly by L. Shepp and A. Shiryaev \cite{MR1233617,MR1348192,MR1846789}.

More elaborate arguments, using optimal stopping, yield:
\begin{equation}
  \label{eq:max2}
  \e\Big[\sup_{s\leq T}|B_s|\Big]\leq \sqrt{2\e T},
\end{equation}
as shown in Dubins and Schwarz \cite{MR89m:60101}.
We also learned from L. Dubins \cite{dubins_private} that
\begin{equation}
  \label{eq:max3}
  \e\Big[\sup_{s\leq T}B_s-\inf_{s\leq T} B_s\Big]\leq \sqrt{3 \e T},
\end{equation}
 and in (\ref{eq:max1}), (\ref{eq:max2}) and (\ref{eq:max3}) the constants are optimal.

To our best knowledge, a bound on the law of the local time similar to (\ref{eq:bound_sup}), has not been derived. We show how to obtain it easily, following the same approach as above, only starting with the martingales given in Corollary \ref{thm:Lharm}.
\begin{prop}
  \label{prop:bound_lt}
Let $m$ be a probability measure on $\re_+$ with $\int_0^\infty x dm(x)<\infty$, and $T$ a stopping time, such that $|B_T|\sim m$ and $(B_{t\land T}:t\ge 0)$ is a uniformly integrable martingale. Denote $\rho_T$ the law of $L_T$. The following bound is true
\begin{equation}
  \label{eq:bound_lt}
\e\Big[\Big(L_T-\overline{\rho}_T^{-1}(p)\Big)^+\Big]\leq \e\Big[\Big(L_{T^m}-\overline{\rho}_{T^m}^{-1}(p^*)\Big)^+\Big],\quad p\in [0,1],
\end{equation}
where $T^m$ is given in Proposition \ref{prop:oy}, the inverses $\overline{\rho}^{-1}_\cdot$ are taken left-continuous and  $p^*=\overline{m}\Big(\overline{m}^{-1}(p)\Big)\ge p$.
\end{prop}
\textbf{Remarks.} For $m$ with no atoms, $p^*\equiv p$. In other words, for $m$ with no atoms, we have $\rho^T\preceq \rho^{T^m}$, where $\rho^{T^m}$ is the image of $m$ through the dual Hardy-Littlewood function $\psi_m$, and ``$\preceq$'' indicates the excess wealth order, defined through
\begin{equation}
  \label{eq:excess}
  \rho_1\preceq \rho_2 \Leftrightarrow \forall p\in [0,1]\;\int_{[\overline{\rho_1}^{-1}(p),\infty)}xd\rho_1(x)\leq \int_{[\overline{\rho_2}^{-1}(p),\infty)}xd\rho_2(x).
\end{equation}
 We point out that the excess wealth order, was introduced recently by Shaked and Shanthikumar in \cite{MR98j:60028} (it is also called the right-spread order, cf. Fernandez-Ponce \textit{et al.} \cite{MR99d:60015}) and studied in Kochar \textit{et al.} \cite{MR2003i:60024}, and the above justifies some further investigation. Since in our case $\e L_T=\e L_{T^m}=\int_0^\infty xdm(x)$, the excess wealth order is equivalent to the TTT and  NBUE orders (see Kochar \textit{et al.} \cite{MR2003i:60024}).\\
{\sl Proof.}
Our proof relies on the martingales given in (\ref{eq:HgenmartL}). Assertion (\ref{eq:bound_lt}) is trivial for $p=1$. It holds also for $p=0$, as it just means that $\e L_T=\e L_{T^m}$, which is true, as both quantities are equal to $\int_0^\infty xdm(x)$. This follows from the fact that $(L_t-|B_t|:t\ge 0)$ is a local martingale and $\e L_{T\land R_n}\nearrow \e L_T$ by monotone convergence, and $\e |B_{T\land R_n}|\to \e |B_T|$ by uniform integrability of $(B_{T\land t}:t\ge 0)$, where $R_n$ is a localizing sequence for $L-|B|$.

Take $p\in (0,1)$, $z=\rho_T^{-1}(p)$ and put $g(x)=\mathbf{1}_{(x> z)}$. Using the optional stopping theorem for the martingale in (\ref{eq:HgenmartL}), we obtain
  \begin{eqnarray}
    \label{eq:pomoc2}
    \e\Big[\Big(L_T-z\Big)^+\Big]&=&\e\Big[|B_T|\mathbf{1}_{(L_T> z)}\Big], \quad\textrm{hence}\nonumber\\
    \e\Big[\Big(L_T-z\Big)^+\Big]&\le&\e\Big[|B_T|\mathbf{1}_{(|B_T|\ge \overline{m}^{-1}(p))}\Big]=\e\Big[|B_{T^m}|\mathbf{1}_{(\varphi_m(L_{T^m})\ge \overline{m}^{-1}(p))}\Big]\nonumber\\
&=&\e\Big[|B_{T^m}|\mathbf{1}_{(L_{T^m}\ge \psi_m(\overline{m}^{-1}(p)))}\Big]\nonumber\\
&\le& \e\Big[|B_{T^m}|\mathbf{1}_{(L_{T^m}> \overline{\rho}^{-1}_{T^m}(p^*))}\Big]=\e\Big[\Big(L_{T^m}-\overline{\rho}^{-1}_{T^m}(p^*)\Big)^+\Big],\nonumber
  \end{eqnarray}
where we have used the fact that $\{x\geq \psi_m(y)\}\subset\{\varphi_m(x)\ge y\}$, which gives
$\p(L_{T^m}\ge \psi_m(\overline{m}^{-1}(p)))\le\p(|B_{T^m}|> \overline{m}^{-1}(p))= p^*$, so that $\psi_m(\overline{m}^{-1}(p))\ge\overline{\rho}_{T^m}^{-1}(p^*)$.
\fdem

\subsection{Penalizations of Brownian motion with a function of its supremum}
\label{sec:penal}
We sketch here yet another instance, where the MM-harmonic functions play a natural role.

Let $f:\re_+\to \re_+$ denote a probability density on $\re_+$, and consider the family of probabilities $(\mathbf{W}^f_t: t\geq 0)$ on $\Omega=C(\re_+,\re)$, where $X_t(\omega)=\omega(t)$, and $\mathcal{F}_s=\sigma(X_u:u\le s)$, $\mathcal{F}_\infty=\bigvee_{s\ge 0} \mathcal{F}_s$, which are defined by:
\begin{equation}
  \label{eq:penal}
  \mathbf{W}^f_t=\frac{f(\overline{X}_t)}{\e_\mathbf{W}\Big[f(\overline{X}_t)\Big]}\cdot \mathbf{W},
\end{equation}
where $\mathbf{W}$ denotes the Wiener measure. Roynette, Vallois and Yor \cite{rvy, rvy2} show that 
\begin{equation}
  \label{eq:penalcnv}
  \mathbf{W}^f_t \xrightarrow[t\to\infty]{(w)}\mathbf{W}_\infty^f,\quad\textrm{i.e.:}\forall s>0,\; \forall \Gamma_s\in \mathcal{F}_s,\; \mathbf{W}^f_t(\Gamma_s)\xrightarrow[t\to\infty]{}\mathbf{W}^f_\infty(\Gamma_s),
\end{equation}
where the probability $\mathbf{W}^f_\infty$ may be described as follows: for $s>0$ and $\Gamma_s\in \mathcal{F}_s$,
\begin{eqnarray}
  \mathbf{W}^f_\infty(\Gamma_s)&=&\e_\mathbf{W}\Big(\mathbf{1}_{\Gamma_s}S^f_s),\quad\textrm{where}\nonumber\\
S^f_s&=&1-F(\overline{X}_s)+f(\overline{X}_s)(\overline{X}_s-X_s)=1-\int_0^s f(\overline{X}_u)dX_u.
\end{eqnarray}
We recognize instantly in the process $S^f$ the max-martingale given by (\ref{eq:Hgenmart}). Another description of $\mathbf{W}^f_\infty$ is that, under this measure the process $X_t$ satisfies:
\begin{equation}
  X_t=X^f_t-\int_0^t\frac{f(\overline{X}_u)du}{1-F(\overline{X}_u)+f(\overline{X}_u)(\overline{X}_u-X_u)},
\end{equation}
where $X^f$ is a $\mathbf{W}^f_\infty$-Brownian motion, and $F(x)=\int_0^x f(u)du$. Naturally, we see the max-martingales (\ref{eq:Hgenmart}) intervene again. Further descriptions of $\mathbf{W}^f_\infty$ are given in  Roynette, Vallois and Yor \cite{rvy2}.

\section{A more general viewpoint: the balayage formula}
\label{sec:balayage}

To end this paper, we propose a slightly more general viewpoint on results mentioned sofar.
In order to present the $(B,\overline{B})$-harmonic functions (\ref{eq:Hgen}), we relied on It\^o's formula. However, it is possible to obtain these functions (and the corresponding martingales) as a consequence of the so-called balayage formula (see, e.g. Revuz and Yor \cite{MR2000h:60050} pp. 260-264 and a series of papers in \cite{MR80f:60003}).

Let $(\Sigma_t:t\ge 0)$ denote a continuous semimartingale, with $\Sigma_0=0$, and define $g_t=\sup\{s\leq t: \Sigma_s=0\}$, $d_t=\inf\{s>t:\Sigma_s=0\}$. Then, the \textit{balayage formula} is: for any locally bounded predictable process $(k_u:u\ge 0)$, one has:
\begin{equation}
  \label{eq:balayage}
  k_{g_t}\Sigma_t=\int_0^tk_{g_s}d\Sigma_s,\quad t\ge 0.
\end{equation}
The intuitive meaning of this formula is that a ``global multiplication'' of $\Sigma$ over its excursions away from $0$ coincides with the stochastic integral of the multiplicator with respect to $(d\Sigma_s)$.
As applications, we give some examples:
\begin{itemize}
\item for $\Sigma_t=\overline{B}_t-B_t$ and $k_u=f(\overline{B}_u)$, $f$ a locally integrable function, (\ref{eq:balayage}) reads $f(\overline{B}_t)(\overline{B}_t-B_t)=\int_0^tf(\overline{B}_s)d(\overline{B}_s-B_s)$, which yields (\ref{eq:Hgenmart});
\item for $\Sigma_t=B_t$ and $k_u=f(L_u)$, $f$ a locally integrable function, $f(L_t)B_t=\int_0^t f(L_s)dB_s$;
\item for $\Sigma_t=|B|_t$ and $k_u=f(L_u)$, $f$ a locally integrable function, we obtain $f(L_t)|B_t|=\int_0^t f(L_s)d|B_s|$. This in turn is equal to $\int_0^t f(L_s)sgn(B_s)dB_s-F(L_t)$ by It\^o-Tanaka's formula, and so we obtain (\ref{eq:HgenmartL}).
\end{itemize}
\smallskip

\section{Closing remarks}
\label{sec:concl}
Max-martingales, or max-harmonic functions, described in (\ref{eq:Hgen}) and (\ref{eq:Hgenmart}), occur in a number of studies of either Brownian motion, or martingales. 
They often lead to simple calculations, and/or formulae, mainly due to the (obvious, but crucial) fact that $d\overline{B}_t$ is carried by $\{t:B_t=\overline{B}_t\}$. 
This has been used again and again by a number of researchers, e.g: Hobson and co-workers, and of course A. Shiryaev and co-workers.
We tried to present in this article several such instances. More generally, this leads to a ``first order stochastic calculus'', as in Section \ref{sec:balayage}, which is quite elementary in comparison with It\^o's second order calculus.

\def\cprime{$'$} \def\polhk#1{\setbox0=\hbox{#1}{\ooalign{\hidewidth
  \lower1.5ex\hbox{`}\hidewidth\crcr\unhbox0}}}

\end{document}